\documentstyle[twoside,psfig,epsf,amsfonts]{amsart}

\begin{document}

\title{MEAN AND VARIANCE ESTIMATION BY KRIGING}

\author{T. SUS{\L}O}

\begin{abstract}
The aim of the paper is to derive the numerical least-squares estimator 
for mean and variance of random variable. 
In order to do so the following questions have to be answered: 
(i) what is the statistical model for the estimation procedure? 
(ii) what are the properties of the estimator, 
like optimality (in which class) or asymptotic properties? 
(iii) how does the estimator work in practice, how compared to 
competing estimators?
\end{abstract}

\maketitle
 
\thispagestyle{empty}

\setcounter{footnote}{0}
\renewcommand{\thefootnote}{\alph{footnote}}

\vspace*{4pt}
\normalsize\baselineskip=13pt  
\section{Introduction}
\noindent
\label{sec:1}
The kriging model based on the stationary random process 
${\mathcal V}=\{V_j;~{\mathbb N}_1 \ni j \supset i=1,\ldots,n\}$ 
with the estimation statistics $\hat{V}_j=\omega^i_j V_i$   
with an unknown constant mean $m$ and variance $\sigma^2$ 
and some correlation function $\rho$ for the asymptotic solution  
$$
\lim_{j \rightarrow \infty} 
E\{[(V_j-m)-(\hat{V}_j-m)]^2\}
= 
\sigma^2 
$$
has the well-known (co-ordinate independent) least-squares disjunction.   
The aim of the paper is to find (on computer) a co-ordinate dependent 
disjunction of kriging model for non-asymptotic solution 
$$
E\{[(V_j-m)-(\hat{V}_j-m)]^2\}
= 
\sigma^2 \ . 
$$

\section{Correlation function estimator}
\label{sec: 2}
\noindent
Following  
\begin{eqnarray}
\gamma(h)=\frac{1}{2}E\{(V_j-V_{j+h})^2\}  
&=&
\frac{1}{2}E\{V_j^2-2V_jV_{j+h}+V_{j+h}^2\} \nonumber \\
&=&
E\{V_j^2\}-E\{V_j V_{j+h}\} \nonumber \\
&=&
E\{V_j^2\}-E^2\{V_j\}-(E\{V_j V_{j+h}\}-E^2\{V_j\}) \nonumber \\
&=&
\sigma^2-C(h) \ge 0 \nonumber
\end{eqnarray}
we get
$$
0 \le C(h)\slash C(0)=C(h)\slash\sigma^2=1-\gamma(h)\slash\sigma^2 = |\rho(h)| 
\ .
$$
Since the correlation function estimator must be non-increasing  
$$
|\hat{\rho}_n(h)|=1-\hat{\gamma}_n(h)\slash\hat{\sigma}^2
$$ 
then only non-decreasing outcomes 
of the experimental semi-variogram 
$$
\hat{\gamma}_n(h)
=\frac{1}{2}\frac{1}{(n-h)}\sum_{j=1}^{n-h}(v_{j}-v_{j+h})^2
\quad
h=0,\ldots,n-1
$$ 
should be taking into consideration for $h \le d$ 
$$
\hat{\sigma}^2
=
\hat{\gamma}_n(d)\ge\hat{\gamma}_n(d-1)\ge\hat{\gamma}_n(d-2)\ge
\ldots\ge\hat{\gamma}_n(1) > \hat{\gamma}_n(0)=0 \ .
$$
Since 
$$
C(h)=E\{V_j V_{j+h}\}-E\{V_j\}E\{V_{j+h}\}
$$
then
$$
\hat{C}_n(h)
=
\frac{1}{n-h}\sum_{j=1}^{n-h}v_jv_{j+h}
-  
\frac{1}{(n-h)^2}
\sum_{j=1}^{n-h}v_j
\sum_{j=1}^{n-h}
v_{j+h} 
$$ 
and the correlation function estimator is also bounded by 
$$
|\hat{\rho}_n(h)|=\hat{C}_n(h)\slash \hat{C}_n(0)
\qquad h=0,\ldots,d \ . 
$$

\section{Central Limit Theorem}
\label{sec: 3}
\noindent
Let us consider a set of 
frozen in time $t$ stationary random processes
$$
{\mathcal V}_t=
\{V_j;~{\mathbb N}_1 \ni j \supset i=1,\ldots,n\}_t;~t \in {\mathbb N}_{n+1}
$$ 
since from the Central Limit Theorem holds 
\begin{eqnarray}
U_{t=n+1}& = & 
\frac{\left[\frac{1}{n}\sum V_j\right]_{t=n+1}-E\{V_j\}_{t=n+1}}
{\sqrt{E\{V_j^2\}_{t=n+1}-E^2\{V_j\}_{t=n+1}}}\sqrt{n} \Rightarrow N(0;1)
\nonumber \\
  & \vdots &  \nonumber \\
U_{t=n+k}& = & 
\frac{\left[\frac{1}{n}\sum V_j\right]_{t=n+k}-E\{V_j\}_{t=n+k}}
{\sqrt{E\{V_j^2\}_{t=n+k}-E^2\{V_j\}_{t=n+k}}}\sqrt{n} \Rightarrow N(0;1)
\nonumber
\end{eqnarray}
an independent set of theirs outcomes 
$$
u_t=\left[\frac{\bar{v}_j-\omega_j^iv_i}
{\sqrt{\omega_j^iv_i^2-(\omega_j^iv_i)^2}}\sqrt{n}\right]_t 
$$
for 
$$
t=n+1,\ldots,n+k
$$
also must follow the Standard Normal distribution.

\begin{figure}[!t] 
\vspace*{13pt}
\centerline{
\psfig{file=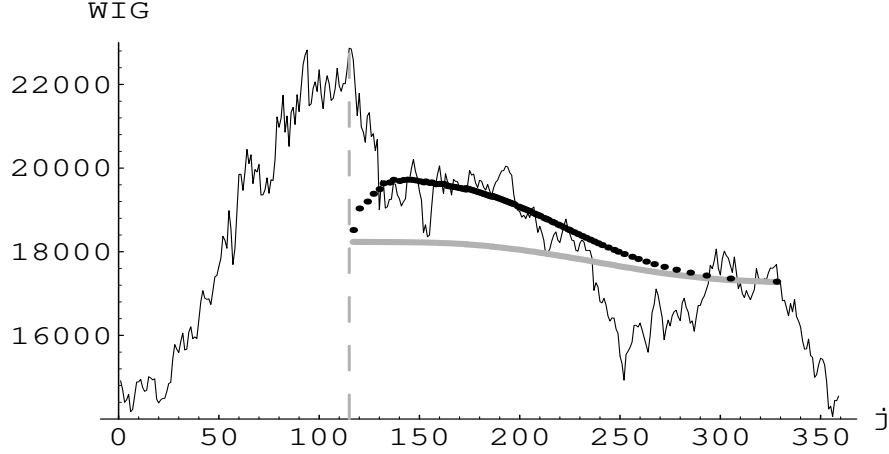,width=12cm,height=6cm}
} 
\vspace*{13pt}
\caption{\label{Fig1} 
Long-lived asymmetric index profile, 
Warsaw's Stock Market Index (WIG) from
11 X 1999 up to 19 III 2001 
(359 close quotes) the numerical least-squares estimator of mean at $j \ge n+1$ 
(black dots) compared for frozen model 
to the asymptotic generalized least-squares
estimator of mean (grey line).  
The dashed vertical line represents $j=n=115$.
} 
\end{figure}

\vspace*{12pt}
\noindent
{\bf Example.} 
Let us consider long-lived asymmetric index profile recorded
by $359$ close quotes of Warsaw's Stock Market Index shown 
in~Fig.~\ref{Fig1}. 
Let the values $v_j=v_{115+1},\ldots,v_{359}$ 
on the latter asymmetric side be supposed to be 
response averages $[\bar{v}_j]_t$ of  
frozen in time $t$ stationary random processes 
based on $v_i=v_1,\ldots,v_{115}$  
$$
{\mathcal V}_t
=\{V_j;~{\mathbb N}_1\ni j \supset i=1,\ldots,115\}_t;~t \in {\mathbb N}_{115+1}
$$    
with the correlation function 
$$
\rho_t(\Delta_{ij})=\left\{
	\begin{array}{ll}
	-1 \cdot {\displaystyle t}^{-1.0135\displaystyle [\Delta_{ij}\slash t]^2},& 
        \qquad \mbox{for}~~\Delta_{ij}>0,\\
	+1, & \qquad \mbox{for}~~\Delta_{ij}= 0,\\
	\end{array}
	\right. 
$$ 
with unknown mean $m_t$ and variance $\sigma^2_t$  
solving on computer the least-squares constraint in one unknown $j$  
for every $t$
$$
\left[ E\{[(V_j-m)-(\hat{V}_j-m)]^2\}
= 
\sigma^2 \  \right]_t \ , 
$$
substituting the solution into  
$$
\hat{u}_t=
\left[\frac{v_j-\omega_j^i v_i}
{\sqrt{\omega_j^i v_i^2-(\omega_j^i v_i)^2}}\sqrt{115}\right]_t 
$$
we get, 
for every $t=115+1,\ldots,115+102$,
the Significance Level for computer K-S test  
$$
\alpha < 0.19 \ .
$$

\end{document}